\documentclass[reqno,12pt]{amsart}
\usepackage{amsmath,amssymb,amsthm,mathtools}
\usepackage{mathrsfs}
\usepackage{bm}
\usepackage{tikz-cd}
\usepackage{booktabs}
\usepackage{xcolor}
\usepackage{hyperref}
\usepackage{float}
\usepackage{orcidlink}
\usepackage{comment}
\usepackage{listings}
\usepackage{xcolor}

\hypersetup{
    colorlinks=true,
    linkcolor=blue,
    citecolor=blue,
    urlcolor=blue}
\usepackage[a4paper,margin=1in]{geometry}
\newtheorem{theorem}{Theorem}[section]
\newtheorem{lemma}[theorem]{Lemma}
\newtheorem{proposition}[theorem]{Proposition}

\theoremstyle{definition}

\allowdisplaybreaks

\title{On the Sum of the Reciprocals of $k$-Generalized Lucas Numbers}

\author{Subhasis Panda} 

\address{Department of Mathematics, School of Advanced Sciences, VIT-AP University, Amravati 522241, Andhra Pradesh, India}

\email{subhasispanda559@gmail.com}

\subjclass[2020]{Primary 11B39, 11B37, 11B83.}

\keywords{Linear recurrences, $k$-generalized Lucas numbers, reciprocal sums}

\begin{document}

\begin{abstract} {For $k\ge2$, let $L^{(k)}$ denote the $k$-generalized Lucas
		sequence. We prove that, for every $n\ge2$, the nearest
		integer to
		$\left(\sum\limits_{m=n}^{\infty}\frac{1}{L_m^{(k)}}\right)^{-1}$
		is $L_n^{(k)}-L_{n-1}^{(k)}$, except for $n=2$ and $n=k+1$.
		In these two cases, the nearest integer is
		$L_n^{(k)}-L_{n-1}^{(k)}-1$ and
		$L_n^{(k)}-L_{n-1}^{(k)}+1$, respectively.
}\end{abstract}

\maketitle

\setcounter{tocdepth}{1}

\section{Introduction}

The study of reciprocal sums of linearly recurrent sequences has
received considerable attention in recent years. The case of Fibonacci numbers was studied by Ohtsuka and
Nakamura \cite{OhtsukaNakamura2008}, who obtained identities
for the integer part of the inverse of infinite reciprocal sums.
Similar identities for Pell numbers were obtained by Zhang
and Wang \cite{ZhangWang2012}. Holliday and Komatsu \cite{HollidayKomatsu2011} extended
these identities to generalized Fibonacci sequences of
order two. Later, Kuhapatanakul
\cite{Kuhapatanakul2013} established identities for alternating
reciprocal sums of generalized Fibonacci numbers.
These results motivate the study of reciprocal sums of
higher-order recurrence sequences.

For an integer $k\ge2$, let $F^{(k)}$ denote the
$k$-generalized Fibonacci sequence defined by
\begin{equation*}
F_n^{(k)}=F_{n-1}^{(k)}+\cdots+F_{n-k}^{(k)} \,\, \text{for all} \,\, n\ge2,	
\end{equation*}
with the initial conditions
$F_{2-k}^{(k)}=\cdots=F_0^{(k)}=0$ and $F_1^{(k)}=1$.
Alahmadi and Luca \cite{AlahmadiLuca2022} proved that for
every $k\ge2$ and $n\ge2$, we have
\begin{equation*}
	\left\langle
	\left(\sum_{m=n}^{\infty}\frac{1}{F_m^{(k)}}\right)^{-1}
	\right\rangle
	=F_n^{(k)}-F_{n-1}^{(k)}.
\end{equation*} 
Their proof relies on the Binet formula of Dresden and Du
\cite{DresdenDu2014} and estimates involving the non-dominant
roots of the characteristic polynomial.

In this paper, we study the corresponding problem for
$k$-generalized Lucas numbers. For an integer $k\geq 2$,
let $L^{(k)}$ denote the $k$-generalized Lucas sequence
defined by the recurrence relation
\begin{equation*}
	L_n^{(k)}=L_{n-1}^{(k)}+\cdots+L_{n-k}^{(k)},
	\qquad n\geq 2,
\end{equation*}
with the initial conditions $L_0^{(k)}=2$, $L_1^{(k)}=1$
and $L_{-j}^{(k)}=0$ for $1\leq j\leq k-2$.
For $k=2$, this sequence reduces to the classical Lucas
sequence. Several arithmetic properties of these numbers
have been studied by Bravo and Luca \cite{BravoLuca2014}
and Bravo et al. \cite{BravoBravoLuca2014}. The $k$-generalized Fibonacci and Lucas sequences are
related by the identity (see \cite[Lemma 1]{BravoLuca2014})
\begin{equation*}
	L_n^{(k)}=2F_{n+1}^{(k)}-F_n^{(k)}.
\end{equation*}
Although the $k$-generalized Fibonacci and Lucas sequences
satisfy the same recurrence relation of order $k$, they
have different initial conditions. The nearest-integer
identity established for $k$-generalized Fibonacci numbers
does not hold for $k$-generalized Lucas numbers at all indices.

In this paper, we determine the nearest integer to
\begin{equation*}
	\left(\sum_{m=n}^{\infty}\frac{1}{L_m^{(k)}}\right)^{-1}
\end{equation*}
for all integers $k\geq 2$ and $n\geq 2$. We prove that
the nearest integer is $L_n^{(k)}-L_{n-1}^{(k)}$,
except when $n=2$ and $n=k+1$, where it is
$L_n^{(k)}-L_{n-1}^{(k)}-1$ and
$L_n^{(k)}-L_{n-1}^{(k)}+1$, respectively. The proof is based on the Binet formula for $k$-generalized
Lucas numbers. We approximate the reciprocal of the infinite
sum by an expression obtained from a geometric series.
The error in this approximation is controlled by the error
term in the Binet formula and its first difference.
To establish suitable upper bounds for these quantities,
we consider two cases $2\leq k\leq 39$ and $k\geq 40$.
These bounds allow us to prove the desired identity for
$n\geq k+2$. The remaining cases are treated separately.

Throughout this paper, we denote by $\langle x\rangle$
the nearest integer to $x$, provided it is unique.
We now state our main result.
\begin{theorem}\label{T1.1}
	\normalfont
	Let $k\geq 2$ and $n\geq 2$ be integers. Then
	\begin{equation*}
		\left\langle
	\left(\sum\limits_{m=n}^{\infty}
	\frac{1}{L_m^{(k)}}\right)^{-1}
	\right\rangle
	=
	\begin{cases}
		L_n^{(k)}-L_{n-1}^{(k)}-1, & \text{if } n=2,\\[2mm]
		L_n^{(k)}-L_{n-1}^{(k)}+1, & \text{if } n=k+1,\\[2mm]
		L_n^{(k)}-L_{n-1}^{(k)}, & \text{otherwise}.
	\end{cases}	
	\end{equation*}
	Moreover, for every $n\geq k+2$, we have
	\begin{equation*}
		\left|
	\left(\sum\limits_{m=n}^{\infty}
	\frac{1}{L_m^{(k)}}\right)^{-1}
	-\left(L_n^{(k)}-L_{n-1}^{(k)}\right)
	\right|< \frac{1}{2}.	
	\end{equation*}
\end{theorem}
The paper is organized as follows. In Section 2, we recall
some known results and prove several preliminary lemmas.
In Section 3, we establish uniform bounds for the error
terms and obtain the required estimate for $n\geq k+2$.
Finally, in Section 4, we prove Theorem~\ref{T1.1}.

\section{Preliminary Results}
Let $\Psi_k(X)=X^k-X^{k-1}-\cdots-X-1$ be the
characteristic polynomial of the $k$-generalized Fibonacci
and Lucas sequences. We also define
\begin{equation*}
	g_k(X)=(X-1)\Psi_k(X)=X^{k+1}-2X^k+1.
\end{equation*}
It is known \cite{Miles1960,Wolfram1998} that $\Psi_k(X)$
is irreducible over $\mathbb{Q}$ and has $k$ distinct roots
$\alpha_1,\alpha_2,\ldots,\alpha_k$, where
$\alpha_1=\alpha>1$ and $|\alpha_i|<1$ for $2\leq i\leq k$.
Moreover, the dominant root $\alpha$ is real and satisfies
$2(1-2^{-k})<\alpha<2$. Since $g_k(z)=0$ for every root
$z$ of $\Psi_k(X)$, we have $z^{k+1}=2z^k-1$.

Let $F^{(k)}$ denote the $k$-generalized Fibonacci sequence
defined by the recurrence relation
$F_n^{(k)}=F_{n-1}^{(k)}+\cdots+F_{n-k}^{(k)}$ for $n\geq 2$,
with the initial conditions
$F_{2-k}^{(k)}=\cdots=F_0^{(k)}=0$ and $F_1^{(k)}=1$.
It follows that $F_n^{(k)}=2^{n-2}$ for $2\leq n\leq k+1$
and $F_{k+2}^{(k)}=2^k-1$. The Binet formula for the terms of $F^{(k)}$, obtained by Dresden and Du \cite{DresdenDu2014}, is given by
\begin{equation*}
	F_n^{(k)}=\sum_{i=1}^k f_k(\alpha_i)\alpha_i^{n-1},
\end{equation*}
where $f_k(z)=\frac{z-1}{2+(k+1)(z-2)}$. We first establish a lower bound for $f_k(\alpha)$,
which will be useful in the proof of our main result.

\begin{lemma}\label{L2.1}
\normalfont For every $k\ge 2$, we have $f_k(\alpha)>\frac12$.
\end{lemma}

\begin{proof}
	Let $t=2-\alpha$. Since $g_k(\alpha)=0$, we have
	$\alpha^{k+1}-2\alpha^k+1=0$, which gives
	$t=\alpha^{-k}$. Since $g_k(1.6)=1-0.4(1.6)^k<0$
	and $g_k(2)=1>0$, we obtain $\alpha>1.6$.
	Moreover, $(k+1)(1.6)^{-k}$ is decreasing for $k\geq 2$.
	Hence
	\begin{equation*}
		(k+1)t=(k+1)\alpha^{-k}
		<(k+1)(1.6)^{-k}
		\leq 3(1.6)^{-2}
		=\frac{75}{64}<\frac65.
	\end{equation*}
	Since $2-(k+1)t>\frac45>0$, we have
	\begin{equation*}
		f_k(\alpha)-\frac12
		=\frac{(k-1)t}{2(2-(k+1)t)}>0.
	\end{equation*}
	Hence $f_k(\alpha)>\frac12$.
\end{proof}

We now recall the Binet formula for $k$-generalized Lucas
numbers and establish a lower bound for the coefficient
of its dominant term.
\begin{lemma}\label{L2.2}
\normalfont Let $k\ge 2$ and $n\ge 2-k$. Then $L_n^{(k)}=h(\alpha)\alpha^{n-1}+\varepsilon_n$, where $h(z)=(2z-1)f_k(z)$,
$|\varepsilon_n|<\frac32$, and $h(\alpha)>\frac{11}{10}$.
\end{lemma}
\begin{proof}
	From Lemma 2 of \cite{BravoLuca2014}, we have
	\begin{equation*}
		L_n^{(k)}
		=\sum_{i=1}^{k}(2\alpha_i-1)
		f_k(\alpha_i)\alpha_i^{n-1}
	\end{equation*}
	and
	\begin{equation*}
		\left|L_n^{(k)}
		-(2\alpha-1)f_k(\alpha)\alpha^{n-1}\right|
		<\frac32 \quad \text{for all } n\geq 2-k.
	\end{equation*}
	Let $h(z)=(2z-1)f_k(z)$. Then we obtain
	\begin{equation*}
		L_n^{(k)}
=h(\alpha)\alpha^{n-1}
+\sum_{i=2}^{k}h(\alpha_i)\alpha_i^{n-1} =h(\alpha)\alpha^{n-1}+\varepsilon_n,		
	\end{equation*}
    where $\varepsilon_n=\displaystyle\sum_{i=2}^{k}
	h(\alpha_i)\alpha_i^{n-1}$ and
	$|\varepsilon_n|<\frac32$ for all $n\geq 2-k$.
	Moreover, since $\alpha>1.6$ and
	$f_k(\alpha)>\frac12$ by Lemma~\ref{L2.1}, we have
	\begin{equation*}
		h(\alpha)=(2\alpha-1)f_k(\alpha)
		>\frac{11}{10}.
	\end{equation*}
\end{proof}

The following lemma follows from 
\cite[Theorem 2.2]{HowardCooper2011}.
\begin{lemma}\label{L2.3}
	\normalfont
	For $k+2\le n\le 2k+2$, we have $F_n^{(k)}=2^{n-2}-(n-k)2^{n-k-3}$.	
\end{lemma}

We use the following lemma of Alahmadi and Luca
\cite[Lemma 3]{AlahmadiLuca2022} to extend the bounds
from a finite range of indices to all larger indices.
\begin{lemma}\label{L2.4}
	\normalfont
	Let $N_0\geq 2-k$ be an integer and let
	$(\delta_n)_{n\geq N_0}$ be a real sequence of the form
	\begin{equation*}
		\delta_n=\sum_{i=2}^{k}c_i\alpha_i^n,
		\qquad n\geq N_0,
	\end{equation*}
	where $c_2,\ldots,c_k$ are complex constants.
	If $|\delta_n|<\lambda$ for $k$ consecutive integers
	$n=n_0,\ldots,n_0+k-1$, where $n_0\geq N_0$
	and $\lambda>0$, then $|\delta_n|<\lambda$
	for all $n\geq n_0$.
\end{lemma}

\section{Bounds for error terms} \label{sec:error-estimates}
In this section, we establish upper bounds for the error
term $\varepsilon_n$ and its first difference.
More precisely, we prove the following proposition.

\begin{proposition}\label{P2.5}
	\normalfont
Let $d_n=\varepsilon_n-\varepsilon_{n-1}$ for $n\ge 1$. For every $k\ge 2$ and $n\ge k+2$, we have
\begin{equation*}
|\varepsilon_n|<\frac14 \,\, \text{and} \,\,
|d_n|<\frac25.	
\end{equation*}
\end{proposition}

We now prove the bounds in Proposition~\ref{P2.5}.
By Lemma~\ref{L2.4}, it is enough to prove both inequalities
for $k+2\leq n\leq 2k+1$. For $2\leq k\leq39$, we verify
both inequalities using SageMath~\cite{SageMath}.
The code used for this verification is given in Appendix~\ref{app:sage}.

To prove the inequalities for $k\geq40$, we first establish
a few lemmas. Let $s_0=2^{-k-1}$ and $s=\alpha^{-k}/2$.
Since $\alpha^k(2-\alpha)=1$, we have $\alpha=2(1-s)$.
Moreover, the bounds $2(1-2^{-k})<\alpha<2$ give
$s_0<s<2^{-k}$. Now, by substituting $\alpha=2(1-s)$ into
$s=\alpha^{-k}/2$, we obtain $s=s_0(1-s)^{-k}$.
For simplicity, we define
\begin{equation*}
	A_n(z)=\frac{(1-2z)(1-z)^{n-1}}{1-(k+1)z}.
\end{equation*}
Using $\alpha=2(1-s)$ and the definition of $f_k$, we have
\begin{align*}
	f_k(\alpha)\alpha^{n-1}
	&=\frac{\alpha-1}{2+(k+1)(\alpha-2)}
	\alpha^{n-1}\\
	&=2^{n-2}\frac{(1-2s)(1-s)^{n-1}}{1-(k+1)s}\\
	&=2^{n-2}A_n(s).
\end{align*}	
From the Binet formula for $F_n^{(k)}$, we have
 $F_n^{(k)}=f_k(\alpha)\alpha^{n-1}+\varepsilon_n^F$, where $\varepsilon_n^F=\displaystyle\sum_{i=2}^{k}
 f_k(\alpha_i)\alpha_i^{n-1}$. Thus, $	\varepsilon_n^F=F_n^{(k)}-2^{n-2}A_n(s)$. We first obtain the Taylor expansion of $A_n(z)$
 about $z=0$ and estimate its remainder term.

\begin{lemma}\label{L3.1}
	\normalfont
For every $k\geq2$ and $n\geq1$, we have
\begin{equation*}
	A_n(z)=1-(n-k)z+c_2(n)z^2+r_n(z),
\end{equation*}
where  \begin{equation*}
	c_2(n)=k^2-kn+k+\frac{n^2-n}{2}-1	
\end{equation*}
and $r_n(z)$ is the remainder term in the Taylor
expansion of $A_n(z)$ about $z=0$.
Moreover, if $k\ge40$ and $k+1\le n\le2k+2$, then $|r_n(z)|\le128(k+1)^3|z|^3$ for every $|z|\le2^{-k}$.
\end{lemma}

\begin{proof}
By taking logarithms of $A_n(z)$ in a neighbourhood of
$z=0$, we obtain
\begin{equation*}
\log A_n(z) =\log(1-2z)+(n-1)\log(1-z) -\log(1-(k+1)z).	
\end{equation*}
From the Taylor expansion of $\log(1-x)$ around
$x=0$, we obtain
\begin{equation*}
\log A_n(z)
=(k-n)z+
\frac{(k+1)^2-(n-1)-4}{2}z^2
+O(z^3).	
\end{equation*}
Now, by taking the exponential of both sides
and using the Taylor expansion of $e^x$, we obtain
\begin{equation*}
	A_n(z) =1+(k-n)z+
	\frac{(k-n)^2+(k+1)^2-(n-1)-4}{2}z^2
	+O(z^3).	
\end{equation*}
Thus, the coefficient of $z^2$ is
\begin{equation*}
		c_2(n)=\frac{(k-n)^2+(k+1)^2-(n-1)-4}{2} =k^2-kn+k+\frac{n^2-n}{2}-1.
\end{equation*}
Therefore, we have $A_n(z)=1-(n-k)z+c_2(n)z^2+r_n(z)$, where $r_n(z)$ is the sum of all terms of degree three and higher.  We now estimate the remainder $r_n(z)$. Let $\rho_0=\frac{1}{2k+2}$. For $k\ge40$,
$n\le2k+2$ and $|z|=\rho_0$, we have
\begin{equation*}
	|1-2z|\le1+\frac{1}{k+1}\le\frac43 \,\, \text{and} \, |1-z|^{n-1}
	\le(1+\rho_0)^{2k+1}<e,
\end{equation*}
where the second inequality follows from $1+x<e^x$ for $x>0$. Also,
$|1-(k+1)z|\ge1-(k+1)|z|=\frac12$. Consequently, we have
\begin{equation*}
	|A_n(z)|
	=\frac{|1-2z|\,|1-z|^{n-1}}
	{|1-(k+1)z|}
	<\frac{8e}{3}<8.
\end{equation*}
Let $A_n(z)=\sum_{j=0}^{\infty}a_jz^j$. Since $A_n(z)$ is analytic in a neighbourhood of
the closed disk $|z|\leq\rho_0$, Cauchy's coefficient
estimate gives
\begin{equation*}
|a_j|\le\frac{8}{\rho_0^j}
=8(2k+2)^j.	
\end{equation*}
Let $x=(2k+2)|z|$. Since $k\ge40$ and
$|z|\le2^{-k}$, we have  $x\le(2k+2)2^{-k}<\frac12$. Therefore, we obtain
\begin{align*}
	|r_n(z)|
	&\le 8\sum_{j=3}^{\infty}x^j
	=\frac{8x^3}{1-x}
	\le16x^3
	=128(k+1)^3|z|^3.
\end{align*}
\end{proof}

We now use Lemma~\ref{L3.1} to derive an explicit approximation
for $\varepsilon_n^F$ and estimate the corresponding error.

\begin{lemma}\label{L3.2}
\normalfont Let $k\ge40$ and define $M_n=-2^{n-2k-5}(n-2k-2)(n-2k+1)$. 	Then, for $k+2\leq n\leq2k+2$, we have $|\varepsilon_n^F-M_n| \le130(k+1)^3 2^{-k}$. Moreover, for $n=k+1$, we have
 \begin{equation*}
 	\left|\varepsilon_n^F-\frac14\right|
 \le2(k+1)^3 2^{-k}.	
 \end{equation*}
\end{lemma}

\begin{proof}
We first establish an inequality that will be used in the proof. Since $k\ge40$, $0<s<2^{-k}$
and $1\le j\le2k$, we have	$js\le2k2^{-k}\le80\cdot2^{-40}<10^{-9}$. Now, using the inequality $-\log(1-s)\le s/(1-s)$, we obtain
\begin{equation*}
0\le-j\log(1-s)
\le\frac{js}{1-s}
<\frac32js<10^{-8}.	
\end{equation*}
Since $e^u-1\le2u$ for $0\le u\le1$,
by taking $u=-j\log(1-s)$, we obtain
\begin{equation}\label{eq1}
	0\le(1-s)^{-j}-1 \le2u \le\frac{2js}{1-s}<3js.
\end{equation}

Let $k+2\le n\le2k+2$. By Lemmas~\ref{L2.3}
and~\ref{L3.1}, we have
\begin{align}
	\varepsilon_n^F
	&=F_n^{(k)}-2^{n-2}A_n(s)\notag\\
	&=2^{n-2}\left[
	(n-k)(s-s_0)-c_2(n)s^2-r_n(s)
	\right].
	\label{eq2}
\end{align}
Since $s_0=s(1-s)^k$, the binomial theorem gives
\begin{equation*}
s-s_0=ks^2-\binom{k}{2}s^3
+\binom{k}{3}s^4-\cdots+(-1)^{k+1}s^{k+1}.
\end{equation*}
Since $ks<1$, the terms after $ks^2$ alternate in sign and decrease in absolute value. By grouping consecutive
terms, we see that their sum lies between $-\binom{k}{2}s^3$ and $0$. Thus, we can write $s-s_0=ks^2-\rho$, where $0\le\rho\le\binom{k}{2}s^3\le\frac{k^2s^3}{2}$. Now, substituting this in equation~\eqref{eq2}, we obtain 
\begin{align}
	\varepsilon_n^F
	&=2^{n-2}\bigl[(n-k)(ks^2-\rho)
	-c_2(n)s^2-r_n(s)\bigr] \notag \\
	&=-2^{n-2}\bigl[(c_2(n)-k(n-k))s^2
	+(n-k)\rho+r_n(s)\bigr]. 	\label{eq3}
\end{align}
Let $P=(n-2k-2)(n-2k+1)$. Since $c_2(n)-k(n-k)=\frac{P}{2}$, the equation~\eqref{eq3} becomes
\begin{equation*}
\varepsilon_n^F
=-2^{n-3}s^2P
-2^{n-2}\bigl[(n-k)\rho+r_n(s)\bigr].	
\end{equation*}
By Lemma~\ref{L3.1}, we have
$|r_n(s)|\le128(k+1)^3s^3$. Hence, we have
\begin{equation*}
\left|\varepsilon_n^F+2^{n-3}s^2P\right| \le2^{n-2}\left[
(n-k)\rho+|r_n(s)|
\right] \le2^{n-2}\left[
\frac{(n-k)k^2}{2}
+128(k+1)^3
\right]s^3.
\end{equation*}
As $n\le2k+2$, we have
\begin{equation*}
\frac{(n-k)k^2}{2}
\le\frac{(k+2)k^2}{2}
\le(k+1)^3	
\end{equation*}
Therefore, $\left|\varepsilon_n^F+2^{n-3}s^2P\right|
\le129(k+1)^3 2^{n-2}s^3$. Since $s<2^{-k}$ and $n\leq2k+2$, we have
$2^{n-2}s^3\leq2^{2k}2^{-3k}=2^{-k}$.
Consequently, we get
\begin{equation}
\left|\varepsilon_n^F+2^{n-3}s^2P\right|
\le129(k+1)^3 2^{-k}.	
\end{equation}
We now replace $s^2$ by $s_0^2$.
Since $s_0=s(1-s)^k$, we have
$s^2=s_0^2(1-s)^{-2k}$.
Applying equation~\eqref{eq1} with $j=2k$, we obtain
$0\leq s^2-s_0^2\leq6ks_0^2s$.
Moreover, for $k+2\leq n\leq2k+2$, we have
$|P|=|(n-2k-2)(n-2k+1)|\leq k^2$.
Now, using $s_0=2^{-k-1}$, we obtain
\begin{equation*}
2^{n-3}|s^2-s_0^2||P|
\le 2^{n-3}\cdot6ks_0^2s\cdot k^2
=6k^3 2^{n-2k-5}s.	
\end{equation*}
Since $n\leq2k+2$, we have $2^{n-2k-5}\leq2^{-3}=1/8$. Combining this with $s<2^{-k}$ gives
\begin{equation*}
6k^3 2^{n-2k-5}s
\le6k^3\cdot\frac18\cdot2^{-k}
=\frac34 k^3 2^{-k}.	
\end{equation*}
Finally, using $\frac34 k^3\le(k+1)^3$, we obtain
\begin{equation*}
2^{n-3}|s^2-s_0^2||P|
\le\frac34 k^3 2^{-k}
\le(k+1)^3 2^{-k}.	
\end{equation*}
Since $M_n=-2^{n-2k-5}P=-2^{n-3}s_0^2P$, by the triangle inequality, we obtain
\begin{align*}
	|\varepsilon_n^F-M_n| & = |\left(\varepsilon_n^F+2^{n-3}s^2P\right)
	+2^{n-3}(s_0^2-s^2)P| \\
	&\le\left|\varepsilon_n^F+2^{n-3}s^2P\right|
	+2^{n-3}|s^2-s_0^2||P|\\
	&\le129(k+1)^3 2^{-k}
	+(k+1)^3 2^{-k}\\
	&=130(k+1)^3 2^{-k}.
\end{align*}

We now consider the case $n=k+1$. In this case,
$F_{k+1}^{(k)}=2^{k-1}$.
By Lemma~\ref{L3.1}, we have
$A_{k+1}(s)=1-s+c_2(k+1)s^2+r_{k+1}(s)$. Hence
\begin{align*}
	\varepsilon_{k+1}^F
	&=F_{k+1}^{(k)}-2^{k-1}A_{k+1}(s)\\
	&=2^{k-1}\bigl[1-A_{k+1}(s)\bigr]\\
	&=2^{k-1}\bigl[s-c_2(k+1)s^2-r_{k+1}(s)\bigr].
\end{align*}
As $s=s_0(1-s)^{-k}$ and $s_0=2^{-k-1}$, we have $2^{k-1}s=\frac14(1-s)^{-k}$. Therefore, we obtain
\begin{equation*}
\varepsilon_{k+1}^F-\frac14
=\frac14\bigl[(1-s)^{-k}-1\bigr]
-2^{k-1}\bigl[c_2(k+1)s^2+r_{k+1}(s)\bigr].	
\end{equation*}
By \eqref{eq1}, we have $0\leq(1-s)^{-k}-1\leq3ks$.
Moreover, from the definition of $c_2(n)$, we have
$c_2(k+1)=\frac{k(k+1)}{2}-1$ and hence
$|c_2(k+1)|\leq(k+1)^2$.
By Lemma~\ref{L3.1}, we have
$|r_{k+1}(s)|\leq128(k+1)^3s^3$. Combining these inequalities, we obtain
\begin{align*}
	\left|\varepsilon_{k+1}^F-\frac14\right| &
	\le\frac14\left|(1-s)^{-k}-1\right|
	+2^{k-1}\left[|c_2(k+1)|s^2
	+|r_{k+1}(s)|\right] \\
	& \le\frac34ks+2^{k-1}
	\left[(k+1)^2s^2+128(k+1)^3s^3\right].	
\end{align*}
Using $s<2^{-k}$, we obtain
\begin{align*}
	\left|\varepsilon_{k+1}^F-\frac14\right|
	&\le\frac34k2^{-k}
	+2^{k-1}\left[(k+1)^2 2^{-2k}
	+128(k+1)^3 2^{-3k}\right]\\
	&=2^{-k}\left[\frac34k+\frac12(k+1)^2
	+64(k+1)^3 2^{-k}\right].
\end{align*}
Since $k\geq40$, we have $64\cdot2^{-k}\leq2^{-34}<1$
and $\frac34k+\frac12(k+1)^2\leq(k+1)^3$.
Therefore,
\begin{equation*}
	\left|\varepsilon_{k+1}^F-\frac14\right|
	\le2(k+1)^3 2^{-k}.
\end{equation*}

\end{proof}

We now use the relation between $\varepsilon_n$ and
$\varepsilon_n^F$ together with the bounds established
in Lemma~\ref{L3.2} to prove the following lemma.

\begin{lemma}\label{L3.4}
\normalfont
Let $k\ge40$ and define $\Delta_k=400(k+1)^3 2^{-k}$. For $x=n-2k$, let
\begin{equation*}
M_n^L=2M_{n+1}-M_n=-2^{x-5}(3x^2+5x-6).
\end{equation*}
Then, for $k+2\leq n\leq2k+1$, we have
$|\varepsilon_n-M_n^L|\leq\Delta_k<10^{-4}$.
Moreover, for $n=k+1$, we have
\begin{equation*}
\left|\varepsilon_{k+1}+\frac14\right|\le\Delta_k.	
\end{equation*}
\end{lemma}

\begin{proof}
Recall that $\varepsilon_n=2\varepsilon_{n+1}^F-\varepsilon_n^F$.
For $k+2\le n\le2k+1$, we have
\begin{align*}
	\varepsilon_n-M_n^L
	&=2\varepsilon_{n+1}^F-\varepsilon_n^F
	-(2M_{n+1}-M_n)\\
	&=2(\varepsilon_{n+1}^F-M_{n+1})
	-(\varepsilon_n^F-M_n).
\end{align*}
By the triangle inequality and Lemma~\ref{L3.2}, we obtain
\begin{align*}
	|\varepsilon_n-M_n^L|
	&\le2|\varepsilon_{n+1}^F-M_{n+1}|
	+|\varepsilon_n^F-M_n|\\
	&\le390(k+1)^3 2^{-k} <\Delta_k.
\end{align*}
Since $\Delta_k$ is decreasing for $k\geq40$, we have
$\Delta_k\leq\Delta_{40}<10^{-4}$.
We now consider the case $n=k+1$. By the identity
$\varepsilon_{k+1}=2\varepsilon_{k+2}^F-\varepsilon_{k+1}^F$
and Lemma~\ref{L3.2}, we obtain
\begin{align}
	\left|\varepsilon_{k+1}+\frac14\right| \notag
	&=\left|2\varepsilon_{k+2}^F-\varepsilon_{k+1}^F+\frac14\right| \notag \\
	&=\left|2(\varepsilon_{k+2}^F-M_{k+2})+2M_{k+2}
	-\left(\varepsilon_{k+1}^F-\frac14\right)\right| \notag\\
	&\le2|\varepsilon_{k+2}^F-M_{k+2}|+2|M_{k+2}|
	+\left|\varepsilon_{k+1}^F-\frac14\right| \notag \\
	&\le262(k+1)^3 2^{-k}+2|M_{k+2}|. \label{eq5}
\end{align}
From the definition of $M_n$, we have
$2|M_{k+2}|=2^{-k-2}k(k-3)\leq k^2 2^{-k}$.
Using equation~\eqref{eq5}, we obtain
\begin{align*}
	\left|\varepsilon_{k+1}+\frac14\right|
	&\le\left[262(k+1)^3+k^2\right]2^{-k}\\
	&<400(k+1)^3 2^{-k}=\Delta_k.
\end{align*}
	
\end{proof}

\noindent\textbf{Proof of Proposition~\ref{P2.5} for $k\ge40$:} 
\begin{proof}
By Lemma~\ref{L3.4} and the triangle inequality, we have
\begin{equation*}
|\varepsilon_n|\leq|\varepsilon_n-M_n^L|+|M_n^L|
<10^{-4}+|M_n^L|	
\end{equation*}
Thus, we first find a uniform upper bound for $|M_n^L|$.
For $k+2\leq n\leq2k+1$, we have $2-k\leq x\leq1$,
where $x=n-2k$. From the definition of $M_n^L$, we have
$|M_n^L|=\frac18,\frac3{16},\frac18,\frac1{32}$
for $x=1,0,-1,-2$, respectively.
For $x\leq-3$, we have
$0<3x^2+5x-6\leq3x^2$, and hence
$|M_n^L|\leq3\cdot2^{x-5}x^2$.
For integers $x\leq-3$, we have
\begin{equation*}
\frac{2^{x-1}(x-1)^2}{2^xx^2}
=\frac12\left(1-\frac1x\right)^2
\le\frac89<1.
\end{equation*}
Thus, $2^xx^2$ is increasing for integers $x\le-3$,
and hence
\begin{equation*}
|M_n^L|\le3\cdot2^{x-5}x^2
\le\frac{27}{256}<\frac3{16}.	
\end{equation*}
Therefore, $|M_n^L|\le3/16$ for $k+2\le n\le2k+1$.
Now, by Lemma~\ref{L3.4}, we obtain
\begin{equation*}
|\varepsilon_n|\le|M_n^L|+\Delta_k
\le \frac3{16}+10^{-4}<\frac14.
\end{equation*}
Next, we find an upper bound for $|d_n|$ for $k+3\leq n\leq2k+1$.
From the definition of $M_n^L$, we have
\begin{align*}
	M_n^L-M_{n-1}^L
	&=-2^{x-6}(3x^2+11x-4)\\
	&=-2^{x-6}(x+4)(3x-1).
\end{align*}
For $x=-4,-3,-2,-1,0,1$, these values are
$0,\frac5{256},\frac7{128},\frac3{32},
\frac1{16},-\frac5{16}$, respectively.
For $x\leq-5$, we have
$0<(x+4)(3x-1)\leq3x^2$.
Since $2^xx^2$ is increasing for integers $x\leq-3$,
we obtain
\begin{equation*}
|M_n^L-M_{n-1}^L|
\le3\cdot2^{x-6}x^2
\le\frac{75}{2048}<\frac5{16}.	
\end{equation*}
By Lemma~\ref{L3.4}, we obtain
\begin{align*} 
	|d_n| & \le \big|(\varepsilon_n-M_n^L)
	+(M_n^L-M_{n-1}^L)
	+(M_{n-1}^L-\varepsilon_{n-1})\big|\\
	&\le|M_n^L-M_{n-1}^L|+2\Delta_k \le\frac5{16}+2 \cdot 10^{-4}<\frac25.
\end{align*}
Finally, we consider the case $n=k+2$. By Lemma~\ref{L3.4}, we obtain
\begin{equation*}
\left|d_{k+2}-\left(M_{k+2}^L+\frac14\right)\right|
\le|\varepsilon_{k+2}-M_{k+2}^L|
+\left|\varepsilon_{k+1}+\frac14\right|
\le2\Delta_k.	
\end{equation*}
From the definition of $M_n^L$, we have
$|M_{k+2}^L|\leq3\cdot2^{-k-3}(k-2)^2
\leq k^2 2^{-k}$.
Since $k^2 2^{-k}$ is decreasing for $k\geq40$, we have
$k^2 2^{-k}\leq40^2 2^{-40}<10^{-8}$.
Therefore, we have
\begin{align*}
	|d_{k+2}|
	&\le\frac14+|M_{k+2}^L|+2\Delta_k\\
	&\le\frac14+k^2 2^{-k}+2\Delta_k\\
	&<\frac14+10^{-8}+2\cdot10^{-4}
	<\frac25.
\end{align*}
Thus, we have $|\varepsilon_n|<\frac14$ and
$|d_n|<\frac25$ for $k$ consecutive integers
$k+2\le n\le2k+1$. Now, by Lemma~\ref{L2.4}, both inequalities hold for all $n\ge k+2$. This completes the proof.
\end{proof}

Recall that $\displaystyle S_n=\sum_{m=n}^{\infty}
\frac{1}{L_m^{(k)}}$ and $D_n=L_n^{(k)}-L_{n-1}^{(k)}$.
We now use the bounds in Proposition~\ref{P2.5}
to show that $D_n$ is the nearest integer to $S_n^{-1}$
for every $k\geq2$ and $n\geq k+2$.
\begin{proposition}\label{P3.5}
\normalfont
Let $k\ge2$ and $n\ge k+2$. Then $\left|S_n^{-1}-D_n\right|<\frac12$.
In particular, $D_n$ is the nearest integer to $S_n^{-1}$.
\end{proposition}

\begin{proof}
By Proposition~\ref{P2.5}, for $k\ge2$ and $n\ge k+2$, we have
$|\varepsilon_m|<1/4$ for every $m\ge n$ and $|d_n|<2/5$. Recall from Lemma~\ref{L2.2} that
$L_m^{(k)}=h(\alpha)\alpha^{m-1}+\varepsilon_m$.
Let $Q_n=h(\alpha)\alpha^{n-1}(1-\alpha^{-1})$. Then
\begin{equation*}
	\sum_{m=n}^{\infty}\frac{1}{h(\alpha)\alpha^{m-1}}
	=\frac{1}{h(\alpha)\alpha^{n-1}}
	\sum_{j=0}^{\infty}\alpha^{-j}
	=\frac{1}{Q_n}.
\end{equation*}
Let $T_n=S_n-\frac{1}{Q_n}$ and $\eta_n=Q_nT_n$. Then $S_n=\frac{1+\eta_n}{Q_n}$ and $S_n^{-1}=\frac{Q_n}{1+\eta_n}$.  Hence
\begin{align}\label{eq6}
	|T_n|
	&=\left|
	\sum_{m=n}^{\infty}
	\left(
	\frac{1}{L_m^{(k)}}-
	\frac{1}{h(\alpha)\alpha^{m-1}}
	\right)
	\right|\nonumber\\
	&=\left|
	\sum_{m=n}^{\infty}
	\frac{-\varepsilon_m}
	{L_m^{(k)}h(\alpha)\alpha^{m-1}}
	\right|.
\end{align}
Let $\kappa_n =\frac{1}{4h(\alpha)\alpha^{n-1}}$. Since $n\ge k+2$ and $k\ge2$, we have $n-1\ge k+1\ge3$. Now, using $\alpha>1.6$ and $h(\alpha)>1.1$, we get
\begin{equation*}
\kappa_n = \frac{1}{4h(\alpha)\alpha^{n-1}}
<\frac{1}{4(1.1)(1.6)^3}<0.056.	
\end{equation*}
For every $m\ge n$, since $\alpha>1$ and
$|\varepsilon_m|<\frac14$, we have
\begin{align}
	L_m^{(k)}
	&=h(\alpha)\alpha^{m-1}+\varepsilon_m
	>h(\alpha)\alpha^{m-1}-\frac14 \notag\\
	&=h(\alpha)\alpha^{m-1}
	\left(1-\frac{1}{4h(\alpha)\alpha^{m-1}}\right)
	\ge h(\alpha)\alpha^{m-1}(1-\kappa_n)>0. \label{eq7}
\end{align}
Applying the lower bound in \eqref{eq7} and the inequality
$|\varepsilon_m|<1/4$ to \eqref{eq6}, we obtain
\begin{align*}
	|T_n|
	&=\left|\sum_{m=n}^{\infty}
	\frac{-\varepsilon_m}{L_m^{(k)}h(\alpha)\alpha^{m-1}}\right|
	\le\sum_{m=n}^{\infty}
	\frac{|\varepsilon_m|}{L_m^{(k)}h(\alpha)\alpha^{m-1}}\\
	&\le\frac{1}{4(1-\kappa_n)h(\alpha)^2}
	\sum_{m=n}^{\infty}\alpha^{-2m+2}
	=\frac{1}{4(1-\kappa_n)h(\alpha)^2
		\alpha^{2n-2}(1-\alpha^{-2})}.
\end{align*}
Since $\eta_n=Q_nT_n$, we have
\begin{align*}
	|\eta_n|
	&=Q_n|T_n|
	\le\frac{h(\alpha)\alpha^{n-1}(1-\alpha^{-1})}
	{4(1-\kappa_n)h(\alpha)^2
		\alpha^{2n-2}(1-\alpha^{-2})}\\
	&=\frac{1}
	{4(1-\kappa_n)h(\alpha)\alpha^{n-1}(1+\alpha^{-1})}
	=\frac{\kappa_n}{(1-\kappa_n)(1+\alpha^{-1})}.
\end{align*}
Since $\alpha<2$ and $\kappa_n<0.056$, we have
$1+\alpha^{-1}>\frac32$ and $1-\kappa_n>0.944$.
Hence, $|\eta_n|<\frac{0.056}{0.944(1.5)}<0.04$.

We now find an upper bound for $|S_n^{-1}-Q_n|$.
Since $S_n^{-1}=\frac{Q_n}{1+\eta_n}$ and
$|\eta_n|<0.04$, we have
\begin{align*}
	|S_n^{-1}-Q_n|
	&=\frac{Q_n|\eta_n|}{|1+\eta_n|}
	\le\frac{Q_n|\eta_n|}{1-|\eta_n|}\\
	&\le\frac{Q_n\kappa_n}
	{(1-\kappa_n)(1+\alpha^{-1})(1-|\eta_n|)}\\
	&=\frac{\alpha-1}
	{4(\alpha+1)(1-\kappa_n)(1-|\eta_n|)}.
\end{align*}
Since $\alpha<2$, we have $\frac{\alpha-1}{\alpha+1}<\frac13$.
It then follows from the bounds $\kappa_n<0.056$ and
$|\eta_n|<0.04$ that
\begin{equation*}
|S_n^{-1}-Q_n|
<\frac{1}{4(3)(0.944)(0.96)}
<0.092.	
\end{equation*}
Finally, by the definitions of $Q_n$, $D_n$ and $d_n$, we have $Q_n=(L_n^{(k)}-\varepsilon_n)
-(L_{n-1}^{(k)}-\varepsilon_{n-1})
=D_n-d_n$. By Proposition~\ref{P2.5},
we obtain
\begin{align*}
	|S_n^{-1}-D_n|
	&\le |S_n^{-1}-Q_n|+|Q_n-D_n|
	=|S_n^{-1}-Q_n|+|d_n|\\
	&<0.092+0.4<\frac12.
\end{align*}
This completes the proof.
\end{proof}

\section{Proof of the main theorem}

\begin{proof}[\normalfont Proof of Theorem~\ref{T1.1}]
By Proposition~\ref{P3.5}, for every $k\geq2$ and $n\geq k+2$, we have $|S_n^{-1}-D_n|<\frac12$. Since $D_n$ is an integer, we obtain $\langle S_n^{-1}\rangle=L_n^{(k)}-L_{n-1}^{(k)}$. Thus, it remains to consider the cases $n=2$, $n=k+1$ and $3\leq n\leq k$.	
	
We first consider the case $n=2$. Since $L_1^{(k)}=1$
and $L_2^{(k)}=3$, we have $D_2=2$.
To prove the first case of Theorem~\ref{T1.1}, it suffices
to show that
\begin{equation*}
	\frac12<S_2^{-1}<\frac32.
\end{equation*}	
Since $L_{m+1}^{(k)}\le 2L_m^{(k)}$ for every $m\ge2$
and $L_2^{(k)}=3$, we have
$L_m^{(k)}\le 3\cdot2^{m-2}$ for every $m\ge2$.
Moreover, $L_{k+1}^{(k)}=3\cdot2^{k-1}-2<3\cdot2^{k-1}$.
Therefore,
\begin{equation*}
	S_2=\sum_{m=2}^{\infty}\frac{1}{L_m^{(k)}}
	>\frac13\sum_{m=2}^{\infty}2^{2-m}
	=\frac23.
\end{equation*}
Consequently, $S_2^{-1}<\frac32$. On the other hand, since $L_{m+1}^{(k)}\ge
L_m^{(k)}+L_{m-1}^{(k)}$ for every $m\ge2$, we have
$L_{m+2}^{(k)}>2L_m^{(k)}$. Since $L_2^{(k)}=3$
and $L_3^{(k)}\ge4$, it follows that
$L_{2j+2}^{(k)}\ge3\cdot2^j$ and
$L_{2j+3}^{(k)}\ge4\cdot2^j$ for every $j\ge0$.
Therefore, we have
\begin{equation*}
S_2=\sum_{j=0}^{\infty}
\left(\frac{1}{L_{2j+2}^{(k)}}
+\frac{1}{L_{2j+3}^{(k)}}\right) \le\left(\frac13+\frac14\right)
\sum_{j=0}^{\infty}2^{-j}
=\frac76<2.	
\end{equation*}
Hence, $S_2^{-1}>\frac12$. Now, combining the above inequalities, we obtain	$\frac12<S_2^{-1}<\frac32$. Thus, $\langle S_2^{-1}\rangle=1=L_2^{(k)}-L_1^{(k)}-1$,
which proves the first case of Theorem~\ref{T1.1}. 

We now consider the case $n=k+1$. Let $P=3\cdot2^{k-1}$. Using Lemma~\ref{L2.3} together with the formula
$F_n^{(k)}=2^{n-2}$ for $2\le n\le k+1$ and the identity
$L_n^{(k)}=2F_{n+1}^{(k)}-F_n^{(k)}$, we obtain
$L_k^{(k)}=P/2$, $L_{k+1}^{(k)}=P-2$ and
$L_{k+2}^{(k)}=2P-5$. Hence $D_{k+1}=P/2-2$
and $D_{k+2}=P-3$. By Proposition~\ref{P3.5}
with $n=k+2$, we have $S_{k+2}^{-1}=P-3+\theta$,
where $|\theta|<1/2$.
Since $S_{k+1}=1/L_{k+1}^{(k)}+S_{k+2}$, we obtain
\begin{equation}\label{eq8}
	S_{k+1}=\frac{1}{P-2}+\frac{1}{P-3+\theta}.
\end{equation}
Since $D_{k+1}=P/2-2$, $P\ge6$ and $|\theta|<1/2$,
it follows from \eqref{eq8} that
\begin{equation*}
	S_{k+1}^{-1}-D_{k+1}-\frac12
	=\frac{P-3+\theta(P-1)}{2(2P-5+\theta)}
	>\frac{P-5}{4(2P-5+\theta)}>0.
\end{equation*}
Hence $S_{k+1}^{-1}-D_{k+1}>1/2$.Similarly, we have
\begin{equation*}
	\frac32-\left(S_{k+1}^{-1}-D_{k+1}\right)
	=\frac{3P-7-\theta(P-3)}{2(2P-5+\theta)}
	>\frac{5P-11}{4(2P-5+\theta)}>0.
\end{equation*}
Hence $S_{k+1}^{-1}-D_{k+1}<3/2$. It follows that $\frac12<S_{k+1}^{-1}-D_{k+1}<\frac32$.
Hence $\langle S_{k+1}^{-1}\rangle
=L_{k+1}^{(k)}-L_k^{(k)}+1$, which proves
the second case of Theorem~\ref{T1.1}. 

Finally, we consider the case $3\leq n\leq k$.
Since $L_m^{(k)}=2F_{m+1}^{(k)}-F_m^{(k)}$ and
$F_m^{(k)}=2^{m-2}$ for $2\le m\le k+1$, we have
$L_m^{(k)}=3\cdot2^{m-2}$ for $2\le m\le k$.
Hence $D_n=3\cdot2^{n-3}$ and
\begin{equation*}
	S_n=\sum_{m=n}^{k}\frac{1}{L_m^{(k)}}+S_{k+1}
	=\frac13\sum_{m=n}^{k}2^{2-m}+S_{k+1}
	=\frac{1}{D_n}-\frac{2}{P}+S_{k+1}.
\end{equation*}
From equation~\eqref{eq8}, we have
\begin{align*}
S_n-\frac{1}{D_n}=\frac{2}{P(P-2)}+\frac{3-\theta}{P(P-3+\theta)}.
\end{align*}
Since $|\theta|<1/2$, both terms on the right-hand side
are positive. Hence $S_n>1/D_n$, which implies
$S_n^{-1}<D_n$. Moreover, since $|\theta|<1/2$, we have $3-\theta<7/2$ 
and $P-3+\theta>P-7/2$. Therefore,
\begin{equation*}
	S_n-\frac{1}{D_n}
	<\frac{2}{P(P-2)}+\frac{7}{P(2P-7)}.
\end{equation*}
Since $3\le n\le k$, we have $k\ge3$, $P\ge12$ and $D_n\le P/4$.
Therefore, we get
\begin{align*}
	D_n-S_n^{-1}
	&=\frac{D_n^2(S_n-1/D_n)}{1+D_n(S_n-1/D_n)}
	<D_n^2\left(S_n-\frac{1}{D_n}\right)\\
	&<\frac{1}{16}\left(\frac{2}{1-2/P}+\frac{7}{2-7/P}\right)
	\le\frac{39}{85}<\frac12.
\end{align*}
Thus, $0<D_n-S_n^{-1}<1/2$, and hence
$\langle S_n^{-1}\rangle=L_n^{(k)}-L_{n-1}^{(k)}$
for $3\le n\le k$. This completes the proof of Theorem~\ref{T1.1}.
\end{proof}

\section*{Declarations}

\subsection*{Funding}
Not applicable.

\subsection*{Disclosure of Interest}

No potential conflict of interest was reported by the author.

\bibliographystyle{amsplain}
\bibliography{mybib}

\end{document}